\documentclass[12pt]{article}
\usepackage{amsmath}
\usepackage{amsfonts}
\newtheorem{theorem}{Theorem}[section]
\newtheorem{lemma}{Lemma}[section] 
\newtheorem{corollary}{Corollary}[section] 
\newtheorem{proposition}{Proposition}[section]

\newcommand{\ts}[1]{\langle #1\rangle}

\title{On locally nilpotent maximal subgroups of the multiplicative group of a division ring }

\author{Bui Xuan Hai 
\\{\small\em Faculty  of Mathematics and Computer Science, University of Natural Science}\\{\small\em VNU - HCM City}\\ {\small\em 227 Nguyen Van Cu str., Dist. 5, Ho Chi Minh City, Vietnam}\\ {\small\em e-mail:bxhai@hcmus.edu.vn}}

\date{25 September 2009}
\rm\em 
\normalsize
\begin{document}
\maketitle
\newcommand{\dpcm}{ \hfill \rule{3mm}{3mm}}
\def\Box{\dpcm}
\newsymbol \lneq  2308   
\newsymbol \lsubset 2320    
\def\xd{\linebreak} 
\begin{abstract} 
Let $D$ be a division ring with the center $F$ and $D^*$ be the multiplicative group of $D$.  In this paper we study locally nilpotent  maximal subgroups of $D^*$. We give some conditions that influence the existence of locally nilpotent maximal subgroups in division ring with infinite center. Also, it is shown that if  $M$ is a locally nilpotent maximal subgroup  that is algebraic over $F$, then either it is the multiplicative group of some maximal subfield of $D$ or it is center by locally finite. If, in addition we assume that $F$ is finite and $M$ is nilpotent,  then the second case cannot occur, i.e. $M$ is the multiplicative group of some maximal subfield of $D$. 
\end{abstract}



\newpage

\section{Introduction}

Let $D$ be a division ring. We denote by $F$ the center  and by $D^*$ the multiplicative group of $D$. In this paper we consider only subgroups of $D^*$ and some time we say about subgroups in a division ring with understanding that they are in fact subgroups under multiplication (subgroups  under addition could not be considered in this paper). It is well-known that if $D$ is commutative, then $D$ may contain no  maximal subgroups. One simple example for that is the field of complex numbers. However, in the case of noncommutativity of $D$, the question on the existence of maximal subgroups in $D$ remains still open and it seems to be difficult. In Section 2  we study this question for the class of division rings that are algebraic over their infinite centers. It is well-known that if a division ring $D$ is algebraic over its center $F$ and $F$ is finite then $D$ is commutative (see Jacobson's Theorem in [4, p.219]). So, we restrict our consideration on division rings, whose centers are infinite. The class of division rings we consider in this section includes the division ring of real quaternions. We show that such division rings don't contain any normal maximal subgroups. As consequence of this fact, we give a series of fields (including the field of complex numbers) whose multiplicative groups don't have any maximal subgroups. We  show that such division rings  have no locally nilpotent maximal subgroups. 
Also, here  we investigate some properties of division rings which influence  the existence of maximal subgroups. 
In Section 3 we give some  characterization of  locally nilpotent maximal subgroups in a division ring $D$ that are algebraic over the  center $F$ of $D$. Here we prove that if $M$ is a locally nilpotent maximal subgroup of $D$ that is algebraic over $F$, then either $M$ is the multiplicative group of some maximal subfield of $D$ or $M$ is center by locally finite. Furthermore,   if $F$ is finite and $M$ is nilpotent,  then the second case cannot occur, i.e. $M$ is the multiplicative group of some maximal subfield of $D$. 

Throughout this paper we use the standard symbols and notation. In particular, if $S\subseteq D$ is a nonempty subset of a division ring $D$ then $C_D(S)$ denotes the centralizer of $S$ in $D$, i. e.
$$C_D(S)=\{x\in D\vert~ xa=ax \mbox{ for all } a\in S\}.$$

If $F\subseteq K$ is a field extension and $a\in K$ is algebraic element over $F$, then we denote by $min(F, a)$ the minimal polynomial of $a$, i.e. the irreducible monic polynomial in $F[X]$ which has $a$ as its root. If every element of $D$ is algebraic over its center $F$, then we say that $D$ is {\em algebraic } over $F$. An element $a\in D$ is said to be {\em radical} over $F$ if there exists some positive integer $n(a)$ depending on $a$ such that $a^{n(a)}\in F$. A subset $S\subseteq D$ is {\em radical} over $F$ if every its element is radical over $F$. 

\section{ A division ring with infinite  center}

In the first we consider a special class of division rings including the division ring of real quaternions. In fact, for a division ring $D$ belonging to this class we  assume  that $D$ contains some algebraic closure of its center. As some interesting application we shall give a series of fields whose multiplicative groups don't have any maximal subgroups. 

\begin{lemma} Let $D$ be  a division ring which  is algebraic over its center $F$ and suppose that $D$ contains an algebraic closure $L$ of $F$. Then, for any element $a\in D$, there exists some element $b\in D^*$ such that $bab^{-1}\in L$. 
\end{lemma}

\noindent
{\bf Proof.} Suppose that $a\in D$ is arbitrary. Denote by $min(F, a)$ the minimal polynomial of $a$ over $F$. Since $L$ is an algebraic closure of $F, min(F, a)$ has some root, say $\omega$ in $L$. Thus $a$ and $\omega$ have the same minimal polynomial over $F$. By Dickson's Theorem  [4, p. 265], there exists some element $b\in D^*$ such that $bab^{-1}=\omega\in L$.\dpcm 

Since the division ring $H$ of real quaternions satisfies the supposition of Lemma 2.1 above, the following result generalizes Theorem 13 in \cite{Ak}.

\begin{theorem}If  $D$ is  a division ring as in Lemma  2.1, then, $D^*$ has no normal maximal subgroups.
\end{theorem}

\noindent
{\bf Proof.} Let $L\subseteq D$ be an algebraic closure of $F$ and suppose that $M$ is a normal maximal subgroup of $D^*$. Then, $D^*/M\simeq \mathbb{Z}_p$ for some prime number $p$. Consider an arbitrary element $a\in D^*$. By Lemma 2.1, there exists some element $b\in D^*$ such that $bab^{-1}\in L$. Since $L$ is an algebraic closure of $F$, the polynomial $f(X)=X^p-bab^{-1}\in L[X]$ has some root $c\in L$; hence 
$$f(c)=c^p-bab^{-1}=0.$$

It follows that $a=b^{-1}c^pb=(b^{-1}cb)^p\in M$. Hence $D^*=M$ that is a contradiction.\dpcm 

The following corollary gives a series of fields (including the field of complex numbers) that contain no maximal subgroups.

\begin{corollary} Every algebraically closed field contains no maximal subgroups.
\end{corollary}

\begin{lemma} Suppose that a finite field extension $F\subset K$ does not have proper intermediate subfields. Then:

(i) Either $K$ is separable over $F$, or

(ii) $K$ is purely inseparable over $F$. Moreover, in this case $char F=p > 0, K$ is radical over $F$ and $[K: F]=p$.
\end{lemma}

\noindent
{\bf Proof.} If  there exists a separable over $F$ element $a\in K\setminus F$, then $K=F(a)$ is separable over $F$. Suppose that every element of $K\setminus F$ is inseparable over $F$. Then $K$ is purely inseparable over $F$. Clearly, in this case we have $char F=p > 0$. Now, consider some element $a\in K\setminus F$. We can find some positive integer $n=n(a)$ depending on $a$ such that $a^{p^n}\in F$ . Suppose that $n$ is a minimal positive integer such that $a^{p^n}\in F$. Setting  $b=a^{p^{n-1}}$, we have $b^p\in F$ and $b\not\in F$. It follows that $K=F(b)$ and $min(F, b)=X^p-b^p$, so $[K: F]=p$. \dpcm

\begin{theorem} If  $D$ is  a division ring as in Lemma  2.1, then $D^*$ contains no locally nilpotent maximal subgroups.
\end{theorem} 

\noindent
{\bf Proof.} If $D$ is a field, then by Corollary 2.1 $D$ contains no maximal subgroups. Now, suppose that $D$ is noncommutative and $M$ is a locally nilpotent maximal subgroup of $D^*$. By [3, Th. 3.2], $M$ is the multiplicative group of some maximal subfield $K$ of $D$. Moreover, $M$ contains $F$. Since $D$ is noncommutative, $K\neq F$ and by [2, Th. 1], there are no proper intermediate subfields of the field extension $F\subset K$ and $Gal(K/F)=\{Id_K\}$. If $a\in K\setminus F$, then $K=F(a)$. Since $a$ is algebraic over $F$, it follows that $[K: F]=[F(a): F] < \infty.$  By Lemma 2.1, there exists some element $b\in D^*$ such that $bab^{-1}\in L$ ($L$ is an  algebraic closure of $F$, lying in $D$). Therefore, $bKb^{-1}=bF(a)b^{-1}\subseteq L$ and $bK^*b^{-1}\subseteq L^*$. Since $K^*=M$ is maximal in $D^*, bK^*b^{-1}$ is maximal in $D^*$. This forces $L^*=bK^*b^{-1}$ and consequently $L=bKb^{-1}$. Hence, one can suppose that $K=L$. In particular, it follows that $F\subset K$ is a normal extension.

By Lemma 2.2, either $K$ is separable over $F$ or $K$ is radical over $F$ and $[K: F]=p=char F > 0$. In the first case, since $K$ is  normal  over $F$, it follows that $F\subset K$ is a Galois extension. Therefore, $|Gal(K/F)|=[K: F]\neq 1$, that is a contradiction. In the last case, for any $u\in D^*$, there exists $v\in D^*$ such that $vuv^{-1}\in L=K$;  hence $u\in v^{-1}Kv$. Since $K$ is radical over $F, u$ is radical over $F$ too. Thus, we have proved that $D$ is radical over $F$. Now, by Kaplansky's Theorem (see [4,  p. 259]), $D$ is commutative, that is again a contradiction.\dpcm

Note that in \cite{Hai-Huynh} it was proved that the division ring of real quaternions does not contain nilpotent maximal subgroups. So, the theorem we have proved strongly generalizes this result. 

Now, suppose that $M$ is a maximal subgroup of $D^*$ and $P$ is the simple subfield of $F$. Denote by $P(Z(M))$ the subfield of $D$ generated by $P\cup Z(M)$. Clearly $P(Z(M))$ is the minimal subfield of $D$ containing $Z(M)$. It was proved in \cite{Hai-Huynh} that $F\subseteq P(Z(M))$. Moreover, if $F$ is infinite then $Z(M)=M\cap F$  iff $P(Z(M))=F$ (see [2, Pro. 1]). Using this fact we can prove the following result.

\begin{proposition} Let $D$ be a  noncommutative division ring with  infinite center $F$. Then the following conditions are equivalent:

(i) $Z(M)=M\cap F$ for every maximal subgroup $M$ of $D^*$.

(ii) $D^*$ contains no maximal subgroups that are multiplicative groups of some division subrings of $D$.
\end{proposition}

\noindent
{\bf Proof.} Suppose that (i) holds and $M$ is a maximal subgroup of $D$ such that $K:=M\cup\{0\}$ is a division subring of $D$. By [2, Pro. 1], $F=P(Z(M))=P(Z(K))=Z(K)$. By [2, Lem. 6], $K^*$ is self-normalized in $D^*$. So $C_D(K)=Z(K)$, hence $C_D(K)=F$. By Double Centralizer Theorem we have $C_D(C_D(K))=K$. It follows that $K=C_D(C_D(K))=C_D(F)=D$, that is a contradiction in view of the maximality of $M=K^*$ in $D^*$. 

Conversely, suppose that $D^*$ contains no maximal subgroups that are multiplicative groups of some division subrings of $D$ and  $M$ is a maximal subgroup of $D$. By setting $K:=M\cup\{0\}$ we have $K\subseteq C_D(Z(M))$. So by maximality of $K^*:=M$ in $D^*$, either $C_D(Z(M))^*=K^*$ or $C_D(Z(M))^*=D^*$. Since by supposition, $K$ is not division subring, we have $C_D(Z(M))^*=D^*$; hence $Z(M)\subseteq F$ and consequently $Z(M)=M\cap F$. \dpcm

\begin{corollary} Let $D$ be a noncommutative division ring that is algebraic over its center $F$. If $Z(M)=M\cap F$ for every maximal subgroup $M$ of $D^*$, then $D$ contains no locally nilpotent maximal subgroups. 
\end{corollary}

\noindent
{\bf Proof.} Suppose that  $Z(M)=M\cap F$ for every maximal subgroup $M$ of $D^*$. If $M$ is a locally nilpotent maximal subgroup of $D^*$, then by [3, Th. 3.2], $M\cup\{0\}$ is the  maximal subfield of $D$ that is a contradiction to the conclusion of Proposition 2.1 above.\dpcm

Note that  in \cite{Hai-Huynh} it was proved that in the division ring of real quaternions for every maximal subgroup $M$ we have $Z(M)=M\cap F$. Hence, in view of Proposition 2.1 and Corollary 2.2 above we  have the following corollary:

\begin{corollary} The division ring $H$ of real quaternions contains no maximal subgroups that are multiplicative groups of some division subrings of $H$. Also, $H$ contains no locally nilpotent maximal subgroups.
\end{corollary}
 
Note that the last assertion of this corollary could be also followed from Theorem 2.2 above.

\section{A division ring with finite center}

Let $D$ be a noncommutative division ring with finite center $F$. In this section we give some characterization of    nilpotent maximal subgroups of $D^*$ that are  algebraic over $F$. 

\begin{lemma} Let $D$ be a noncommutative division ring with center $F$ and suppose that $M$ is a locally nilpotent maximal subgroup of $D^*$ that is algebraic over $F$. Then,  one of the following cases occurs:

(i) Either $F^*\subseteq M$ and there exists a maximal subfield $K$ of $D$ such that  $M=K^*$ or,

(ii) $M$ is center by locally finite.
\end{lemma}

 \noindent
{\bf Proof.} Since $M$ is maximal in $D^*$, either $M=F(M)^*$ or $F(M)=D$. If $M=F(M)^*$, then $F^*\subseteq M$ and by [3, Th. 2.2] $K=F(M)$ is the maximal subfield of $D$. 

Now, suppose that $F(M)=D$. Since $M$ is algebraic over $F$, we have $D=F(M)=F[M]$; so $M$ is absolutely irreducble. By [5, Th. 5.7.11, p. 215], $M$  is center by locally finite. \dpcm

Now, we are ready to prove the following result for a division ring with finite center.

\begin{theorem} Let $D$ be a noncommutative division ring with center $F$ and suppose that $M$ is a nilpotent maximal subgroup of $D^*$ that is algebraic over $F$. If $F$ is finite, then $M$ is the multiplicative group of some maximal subfield of $D$.
\end{theorem}

\noindent
{\bf Proof.}
Suppose that $F$ is finite. In view of Lemma  3.1, it suffices to show that the case (ii) cannot occur. Thus, suppose that $F(M)=D$. Since $M$ is maximal in $D^*$, we can show that $F^*\subseteq M$. If not, we have $F^*M=D^*$, so $D'=M'\subseteq M$. Then $D^*$ is solvable and by Hua's Theorem (see, for example [4, p. 223]) it follows that $D$ is commutative, that is a contradiction. Therefore $F^*\subseteq M$, so $F^*=Z(M)$. By Lemma 3.1 (ii)  $M/F^*$ is locally finite. Consider arbitrary elements $x, y\in M$. Then, the subgroup $\ts{xF^*, yF^*}$ of $M/F^*$ generated by $xF^*$ and $yF^*$ is finite. Suppose that $g$ is the restriction of the natural homomorphism $M\longrightarrow M/F^*$ on the subgroup $\ts{x, y}$. Then, we have $Ker g=\ts{x, y}\cap F^*$ and $Im g=\ts{xF^*, yF^*}$. Since $F^*$ and $\ts{xF^*, yF^*}$ are both finite, it follows that $\ts{x, y}$ is finite. Therefore $\ts{x, y}$ is cyclic and in particular, $x, y$ commute with each other. So, $M$ is abelian and consequently, $D=F(M)$ is commutative,  that is a contradiction.\dpcm

\end{document}